\documentclass[11pt]{article}

\usepackage{amsfonts,amsmath,amssymb}
\usepackage{graphicx}
\usepackage{color}
 \usepackage{fullpage}

 \def \qed {\nopagebreak{\hfill$\Box$\medskip}}

 \newtheorem{theorem}{Theorem}
 \newtheorem{proposition}[theorem]{Proposition}
 \newtheorem{corollary}[theorem]{Corollary}
 \newtheorem{conjecture}[theorem]{Conjecture}
 \newtheorem{lemma}[theorem]{Lemma}

 \newtheorem{remark}{Remark}
 \newtheorem{defi}[theorem]{Definition}

 \newcommand{\bt}{\begin{theorem}}
 \newcommand{\et}{\end{theorem}}
 \newcommand{\bl}{\begin{lemma}}
 \newcommand{\el}{\end{lemma}}
 \newcommand{\bp}{\begin{proposition}}
 \newcommand{\ep}{\end{proposition}}
 \newcommand{\bcor}{\begin{corollary}}
 \newcommand{\ecor}{\end{corollary}}
 \newcommand{\br}{\begin{remark}\rm}
 \newcommand{\er}{\end{remark}}
 \newcommand{\bcon}{\begin{conjecture}}
 \newcommand{\econ}{\end{conjecture}}
 \newcommand{\bd}{\begin{defi}}
 \newcommand{\ed}{\end{defi}}
 \newcommand{\ben}{\begin{enumerate}}
 \newcommand{\een}{\end{enumerate}}

 \newcommand{\beq}{\begin{equation}}
 \newcommand{\eeq}{\end{equation}}

 \newcommand{\Bi}{{\cal B}}

 \newcommand{\eps}{\epsilon}

 \newcommand{\GG}{\mbox{${\mathcal G}$}}
 \newcommand{\FF}{\mbox{${\mathcal F}$}}

 \renewcommand{\SS}{\mbox{${\mathcal S}$}}

 \newcommand{\Nbold}{\mbox{${\mathbb N}$}}

 \newcommand{\Zbold}{\mbox{${\mathbb Z}$}}

 \newcommand{\Ebb}{\mbox{${\mathbb E}$}}
 \newcommand{\Pbb}{\mbox{${\mathbb P}$}}

 \newcommand{\bE}{{\bf E}}
 \newcommand{\bP}{{\bf P}}

 \newcommand{\bone}{{\bf 1}}

 \newcommand{\Z}{{\mathbb Z}}

 \newcommand{\T}{{\mathbb T}}
 \newcommand{\W}{{{\mathbb N}_0}}
 \renewcommand{\P}{P}

 \renewcommand{\ln}{\log}

\begin{document}

 \title{{\bf Random Walks in I.I.D. Random Environment on Cayley Trees}}
 \author{{\bf Siva Athreya}$^{\,1}$, \and {\bf Antar Bandyopadhyay}$^{2,3}$, \and {\bf Amites Dasgupta}$^{3}$}
 \date{\today}

 \maketitle

 \footnotetext[1]{Indian Statistical Institute, 8th Mile Mysore Road, Bangalore 560059, India. Email: athreya@isibang.ac.in}

 \footnotetext[2]{Indian Statistical Institute, 7 S. J. S. Sansanwal Marg, New Delhi  110016, India. Email: antar@isid.ac.in }

 \footnotetext[3]{Indian Statistical Institute, 203 Barrackpore Trunk Road, Kolkata 700 108, India. Email: amites@isical.ac.in}

 \medskip

 \begin{abstract}
 We consider the random walk in an \emph{i.i.d.} random environment on the infinite $d$-regular tree for
$d \geq 3$. We consider the tree as a Cayley graph of free product of finitely many copies of $\Zbold$
and $\Zbold_2$ and define the i.i.d. environment as invariant under the action of this group. 
Under a mild non-degeneracy assumption we show that the walk is always transient.

 \vspace{0.15in}


 \vspace{0.15in}

 \end{abstract}

 \section{Introduction}
 \label{Sec:Intro}

In this short note we consider a random walk in random environment (RWRE)
model on a regular tree with degree $d \geq 3$, where the environment at the vertices 
are \emph{independent} and are also ``\emph{identically distributed}'' (i.i.d.). 
We make this notion of \emph{i.i.d. environment}
rigorous by first defining a translation invariant model on a group $G$ which is a free
product of finitely many groups, $G_1, G_2, \cdots, G_k$ and
$H_1, H_2, \cdots, H_r$, where 
each $G_i \cong \Zbold$ and each $H_j \cong \Zbold_2$ with $d = 2k + r$.  
Observing the fact that the Cayley graph of $G$ is a regular tree with degree $d$,
we transfer back the model on the $d$-regular tree we started with. 
We prove that under a mild non-degeneracy assumption such a walk is always
transient.

\subsection{Basic Setup}
\label{SubSec:Setup}

{\bf Cayley graph:} Let $G$ be a group defined above, that is, $G$ is a free product of $k + r \geq 2$
groups, namely
$G_1, G_2, \ldots, G_k$ with $k \geq 0$ and
$H_1, H_2, \ldots, H_r$ with $r \geq 0$, where 
each $G_i \cong \Zbold$ and each $H_j \cong \Zbold_2$ and $d = 2k + r \geq 3$. 
Suppose $G_i = \langle a_i \rangle$ for $1 \leq i \leq k$ and 
$H_j = \langle b_j \rangle$ where $b_j^2 = e$ for $1 \leq j \leq r$.
Here by $\langle a \rangle$ we mean the group generated by a single element $a$.  
Let 
$S := \left\{a_1, a_2, \ldots, a_k\right\} \cup \left\{a_1^{-1}, a_2^{-1}, \ldots, a_k^{-1}\right\}
      \cup \left\{ b_1, b_2, \ldots, b_r \right\}$
be a generating set for $G$. We note that $S$ is a symmetric set, that is, $s \in S \iff s^{-1} \in S$. 
%
%

We now define a graph $\bar{G}$ with vertex set $G$ and edge set 
$E := \left\{ \left\{x,y\right\} \,\Big\vert\, yx^{-1} \in S \,\right\}$. Such a 
graph $\bar{G}$ is called a \emph{(left) Cayley Graph} of $G$ with respect to the generating set $S$. 
Since $G$ is a free product of groups which are isomorphic to either $\Zbold$ or $\Zbold_2$, 
it is easy to see
that $\bar{G}$ is a graph with no cycles and is regular with degree $d$, thus it is isomorphic to 
the $d$-regular infinite tree which we will denote by $\T_d$. 
We will abuse the terminology a bit and will write $\T_d$ for 
the Cayley graph of $G$. 
We will consider the identity element $e$ of $G$ as the root of $\T_d$.
We will write $N\left(x\right)$ for the set of all neighbors of a vertex $x \in \T_d$. Notationally, 
$N\left(x\right) = \left\{ y \in G \,\Big\vert\,  yx^{-1} \in S \,\right\}$. Observe that from definition 
$N\left(e\right) = S$.
For $x \in G$, define the mapping $\theta_x : G \rightarrow G$ by $\theta_x\left(y\right) = yx$, 
then $\theta_x$ is an automorphism of $\T_d$. We will call
$\theta_x$ the \emph{translation by} $x$. For a vertex $x \in \T_d$ and $x \neq e$, 
we denote by $\left\vert x \right\vert$, 
the length of the unique path from the root $e$ to $x$ and $\left\vert e \right\vert = 0$. Further,
if $x \in \T_d$ and $x \neq e$ then we define $\overleftarrow{x}$ as the 
\emph{parent} of $x$, that is, the penultimate vertex on the unique path from
$e$ to $x$.

\medskip

{\bf Random Environment:} Let $\SS := \SS_{e}$ be a  
collection of probability
measures  on the $d$ elements of $N\left(e\right) = S$. 
To simplify the presentation and avoid various mesurability issues,
we assume that $\SS$ is a Polish space (including the possibilities
that $\SS$ is finite or countably infinite).
For each $x \in \T_d$,  $\SS_{x}$ is the the push-forward of the space
$\SS$  under the translation $\theta_x$, that is, 
$\SS_x := \SS \circ \theta_x^{-1}$. 
Note that the probabilities on $\SS_x$ have support on
$N\left(x\right)$. That is to say, an element
$\omega(x,\cdot)$ of $\SS_{x}$, is a probability measure
satisfying
\[
\begin{array}{rcl}
\omega\left(x, y\right) \geq 0 \,\,\,\, \forall \,\,\, 
y \in \T_d
& 
\mbox{and}
& 
\mathop{\sum}\limits_{y \in N\left(x\right)} \omega\left(x, y\right) = 1.
\end{array}
\]

Let ${\cal B}_{{\mathcal S}_x}$ denote the Borel $\sigma$-algebra on $\SS_x$. The 
\emph{environment space} is defined as the measurable space $\left(\Omega, \FF\right)$ where
\begin{equation}
\begin{array}{cc}
\Omega := \mathop{\prod}\limits_{x \in \T_d} \SS_{x},
& 
\FF := \mathop{\bigotimes}\limits_{x \in \T_d} {\cal B}_{{\mathcal S}_x}.
\end{array}
\label{Equ:Omega-FF}
\end{equation}
An element $\omega \in \Omega$ will be written as
$\left\{ \omega\left(x, \cdot\right) \,\Big\vert\, x \in \T^d \,\right\}$.
An environment distribution is a probability $\P$ on $\left(\Omega, \FF\right)$. 
We will denote by $E$ the expectation taken with respect to the probability measure $\P$.

\medskip

{\bf Random Walk:} Given an environment 
$\omega \in \Omega$, a random walk $\left(X_n\right)_{n \geq 0}$
is a time homogeneous Markov chain
taking values in $\T_d$ 
with transition probabilities
\[
\bP_{\omega} \left( X_{n+1} = y \,\Big\vert\, X_n = x \right) = \omega\left(x, y\right). 
\]
Let $\Nbold_0 := \Nbold \cup \left\{ 0 \right\}$. 
For each $\omega \in \Omega$, we denote by
$\bP_{\omega}^{x}$ the law induced by $\left(X_n\right)_{n \geq 0}$ on
$\left( \left(\T_d\right)^{\Nbold_0}, \GG \right)$,
where $\GG$ is the $\sigma$-algebra generated by the cylinder sets, 
such that
\begin{equation}
\bP_{\omega}^{x}\left( X_0 = x \right) = 1. 
\label{Equ:RW-Intitial-Prob}
\end{equation}
The probability measure $\bP_{\omega}^{x}$ is called the \emph{quenched law}
of the random walk $\left(X_n\right)_{n \geq 0}$, starting at $x$.
We will use the notation $\bE_\omega^x$ for the expectation under the quenched measure $\bP_\omega^x$. 

Following Zeitouni \cite{Zei04},
we note that for every $B \in \GG$, the function
\[
\omega \mapsto \bP_{\omega}^{x}\left(B\right)
\]
is $\FF$-measurable. Hence, we may define the measure
$\Pbb^{x}$ on 
$\left( \Omega \times \left(\T_d\right)^{\Nbold_0},
\FF \otimes \GG \right)$ by the relation
\[
\Pbb^{x}\left( A \times B \right)
= 
\int_A \! \bP_{\omega}^{x}\left(B\right) \P\left(d\omega\right), 
\,\,\,\, \forall \,\,\, A \in \FF, \, B \in \GG. 
\]
With a slight abuse of notation, we also denote the marginal of 
$\Pbb^{x}$ on $\left(\T_d\right)^{{\mathbb N}_0}$
by $\Pbb^{x}$, whenever no confusion occurs. This probability distribution
is called
the \emph{annealed law} of the random walk $\left(X_n\right)_{n \geq 0}$, 
starting at $x$. 
We will use the notation $\Ebb^x$ for the expectation under the annealed measure $\Pbb^x$

\subsection{Main Results}
\label{SubSec:Results}

Throughout this paper we will assume that the following hold:

\begin{itemize}
\item[(A1)] $\P$ is a product measure on $\left(\Omega, \FF\right)$ with ``\emph{identical}'' marginals, that is, 
            under $\P$ the random probability laws $\left\{ \omega\left(x, \cdot\right) \,\Big\vert\, x \in \T^d \,\right\}$
            are independent and ``identically'' distributed in the sense that
            \beq \label{si} \P \circ \theta_{x}^{-1} = \P, \eeq
            for all $x \in G$. 
            
\item[(A2)] For all $1 \leq i \leq d$,
            \beq \label{Equ:UEllip-} E\left[\left\vert \log \omega\left(e,s_{i}\right) \right\vert \right] < \infty. \eeq
            It is worth noting that under this assumption $\omega\left(x,y\right) > 0$ almost surely (a.s.)
            with respect to the measure $\P$ for all $x \in \T_d$ and $y \in N\left(x\right)$. 

\end{itemize}

The following is our main result.
\bt 
\label{main}
Under assumptions (A1) and (A2) the random walk $\left(X_{n}\right)_{n \geq 0}$ is  transient $\Pbb^{e}$-a.s.,
that is, 
$\Pbb^{e} \left( \mathop{\lim}\limits_{n \rightarrow \infty} \left\vert X_n \right\vert = \infty \right) = 1$
\et

An immediate question that arises is whether the above walk has a
speed which may be zero. The following result provides a partial answer to
this question with (A2) replaced by the usual \emph{uniform ellipticity} condition.      
\begin{itemize}
\item[(A3)] There exists $\eps > 0$ such that 
            \beq 
            \label{Equ:UEllip} 
            \P\left( \omega\left(e, s_i\right) > \eps \,\, \forall \,\, 1 \leq i \leq d \right) = 1.
            \eeq
\end{itemize}
\bt
\label{Thm:Speed}
Under assumptions (A1) and (A3) with $\eps > \frac{1}{2\left(d-1\right)}$ we have $\Pbb^{e}$-a.s.
\beq
\liminf_{n \rightarrow \infty} \frac{\left\vert X_n \right\vert}{n} > 0.
\label{Equ:Speed}
\eeq
\et
Note that the condition $\eps > \frac{1}{2\left(d-1\right)}$ is compatible with the ellipticity 
condition \eqref{Equ:UEllip} as $d \geq 3$. 

\subsection{Remarks}
\label{Sec:Final}

Random walk in Random Environment (RWRE) model on the one dimensional
integer lattice $\Z$ was first introduced by Solomon in \cite{Sol75}
where he gave explicit criteria for the recurrence and transience of
the walk for i.i.d.
environment.  Since then a large variety of results have
been discovered for RWRE in $\Z^d$, yet there are many challenging
problems which are still left open (see \cite{Bog06} for a non-technical survey and 
\cite{Zei04, Szn04} for more technical details).  

Perhaps the earliest known results for RWRE on trees are by Lyons and Pemantle
\cite{LyPe92}. In that paper they consider a model on rooted
trees known as \emph{random conductance model}. In
that model, the random conductances along each path from vertices to
the root are assumed to be independent and identically
distributed. The random walk is then shown to be recurrent or
transient depending on how large the value of the average
conductance is. 

In our set up, the assumption (A1) essentially says that the
random transition laws $\left\{ \omega\left(x, \cdot\right)
\,\Big\vert\, x \in \T^d \,\right\}$ are \emph{independent and
  identically distributed (i.i.d.)}. On $\T_d$ we introduced the group
structure to define \emph{identically distributed} and we made the
probability law $\P$ invariant under \emph{translations} by the
group elements.  Hence the RWRE model in this article is different
from the \emph{random conductance model} discussed above. It
is interesting to note that the only example where the two models
agree is the deterministic environment of the \emph{simple symmetric walk} on $\T_d$. 

Perhaps the model closest to ours was introduced by Rozikov \cite{Rozi01} where the author 
considered the case with $k=0$ and $r=d \geq 3$, that is, the group $G$ is a free
product for $d$ copies of $\Zbold_2$. Our model is slightly more general from
this perspective, but more importantly in \cite{Rozi01} to prove transience, it was assumed that
\begin{equation}
E\left[ \left\vert \log \frac{\omega\left(x, sx\right)}{\omega \left(x, s'x\right)} \right\vert \right] < \infty
\label{Equ:Rozikov-1}
\end{equation}
and
\begin{equation}
E\left[ \log \frac{\omega\left(x, sx\right)}{\omega \left(x, s'x\right)} \right] \neq 0, 
\label{Equ:Rozikov-2}
\end{equation}
for every $x \in \T_d$ and for two different elements $s, s' \in S$
(see the assumption made in the beginning of Section 7 and Theorem 5 of \cite{Rozi01}). 
The first assumption~\eqref{Equ:Rozikov-1} made in \cite{Rozi01} is more general than
our assumption (A2). However, the second assumption made in \cite{Rozi01}, 
namely equation~\eqref{Equ:Rozikov-2}, may not be satisfied by certain environments 
(be it random or non-random) which
are covered by our formulation, 
for example, the condition~\eqref{Equ:Rozikov-2} is not satisfied by the 
\emph{simple symmetric random walk} on $\T_d$. 
So neither
our model is a subclass of the models studied by Rozikov \cite{Rozi01}, nor our model covers all cases discussed 
in there. 
So we believe our work is an important addition to the
earlier work of Rozikov \cite{Rozi01} and none makes the other redundant. We would also like to point out 
that the techniques used in our work are entirely different from that of \cite{Rozi01}.

There have also been several other contributions on random trees,
particularly on random walk on Galton-Watson trees 
\cite{Ly90, LyPePe95, LyPePe96, DemGanPeZei02, PeZei08}.  It is worth pointing
out here that a random walk on a Galton-Watson tree \cite{Ly90}
satisfies the assumption (A1) and so does a random walk on a
multi-type Galton-Watson tree \cite{DemSun12}.

Our last result (Theorem 2) is certainly far from
satisfactory. We strongly believe that under the assumptions (A1) and
(A3) the sequence of random variables $\left( \frac{\left\vert X_n
  \right\vert}{n} \right)_{n \geq 0}$ has a $\Pbb^e$-almost sure limit
which is non-random and strictly positive.  A similar conclusion has
been derived for the special case of random walk on Galton-Watson
trees \cite{LyPePe95}. This and the central limit theorem for such
walks will be studied in future work.


\section{Proofs of the Main Results}
\label{Sec:Proofs}

\subsection{Proof of Theorem~\ref{main}}
\label{SubSec:Proof-Main}
Given an environment $\omega \in \Omega$, and a vertex $\sigma \in \T_d$ which is not the root, 
we define the conductance of the edge $\left\{ \sigma, \overleftarrow{\sigma} \right\}$ as 
\begin{equation}
C\left( \sigma, \overleftarrow{\sigma} \right)
= \omega (e, x_1) 
\left( 
\prod_{k=1}^{\left\vert \sigma \right\vert-1} \frac{\omega(x_{k},x_{k+1})}{\omega(x_k,x_{k-1})} 
\right),
\end{equation}
where $e = x_0, x_1, x_2, \cdots, x_{\left\vert \sigma \right\vert - 1} = \overleftarrow{\sigma}, 
x_{\left\vert \sigma \right\vert} = \sigma$ is the unique path from the root $e$ to the vertex $\sigma$. 
Further we define $\Phi\left(\sigma\right) := C\left( \sigma, \overleftarrow{\sigma} \right)^{-1}$. 
%
Suppose $\sigma = \alpha_n \alpha_{n-1} \cdots \alpha_1$ where 
$\alpha_{i} \in S$ and $\alpha_{i+1} \neq \alpha_i^{-1}$ and $n = \left\vert \sigma \right\vert$.
More generally $x_k = \alpha_k \alpha_{k-1} \cdots \alpha_1$ for $1 \leq k \leq n$. Note that
$x_k = \alpha_k x_{k-1}$ for all $1 \leq k \leq n$. 
We can now rewrite $\Phi\left(\sigma\right)$ as
\begin{equation}
\Phi\left(\sigma\right) 
= \frac{1}{\omega(e, x_1)} \prod_{k=1}^{n-1} \frac{\omega(x_{k},x_{k-1})}{\omega(x_{k},x_{k+1})}
= \left(
  \prod_{k=1}^{n-1} \frac{\omega_{k}\left(\alpha_{k}^{-1}\right)}{\omega_{k-1}\left(\alpha_{k}\right)} 
  \right)
  \frac{1}{\omega_{n-1}\left( \alpha_n\right)}, 
\end{equation}
where we write $\omega_{k} \left(s\right) := \omega\left(x_{k}, s x_{k}\right)$ for any $s \in S$. 

We will now show that there is a (non-random) sequence of positive real numbers 
$\left(\beta_n\right)_{n \geq 1}$ such that
$\displaystyle{\sum_{n=1}^{\infty} \beta^n < \infty}$ and $P$-a.s.
\begin{equation}
\label{Equ:Main-Criterion}
\lim_{n \rightarrow \infty} \sum_{\sigma_{n} \in \T^{n}_{d}}  \beta_{n} \left(\Phi\left(\sigma_{n}\right)\right)^{-1} = \infty,
\end{equation}
where $\T^{n}_{d} := \left\{ x \in \T_d \,\Big\vert\, \left\vert x \right\vert = n \,\right\}$.
Then by Corollary 4.2 in \cite{Ly90}, the random walk has to be transient.
For this we will show that $P$-a.s., there is a subset of vertices of $\T_d^n$ with size 
$O\left(\left(d-1\right)^{n-1}\right)$ such that the $\Phi$-value of these vertices are strictly smaller
than $\left(d-1\right)^{\frac{n}{2}}$. 

Let $\Bi_\W $  denote the product
$\sigma$-algebra on $S^\W$,  and  $\mu $ be a probability measure on $\left(S^{\W}, \Bi_{\W}\right)$ 
such that $\left(Y_n\right)_{n \geq 0} \in S^{\W}$
forms a Markov chain on $S$ with  
\begin{equation}
\mu\left( Y_{n} = s \,\Big\vert\, Y_{n-1} = t \,\right) = \frac{1}{d-1},  \,\, s, t \in S \mbox{ with } s \neq t^{-1}. 
\label{Equ:MC-Transition-Y}
\end{equation}
It is easy to see that the chain $\left(Y_{n}\right)_{n \geq 0}$ 
is an aperiodic, irreducible and finite state Markov chain and its stationary distribution is 
the uniform distribution on $S$. We shall assume that $Y_0$ is uniformly distributed on $S$.
Thus each $Y_n$ is also uniform on $S$. 

Let $\eta_{n} = Y_n Y_{n-1} \cdots Y_1$. From equation~\eqref{Equ:MC-Transition-Y} 
it follows that 
$\eta_{n}$ is uniformly 
distributed on the set of vertices 
$\T^{n}_{d}$.  Now 
\begin{eqnarray}
\frac{1}{n} \ln \Phi\left(\eta_{n}\right) 
& = &  o\left(1\right) + 
       \frac{1}{n}\sum_{k=1}^{n-1} 
                  \left( \ln \omega_{k}\left(Y_{k}^{-1}\right) - \ln \omega_{k-1}\left(Y_{k}\right) \right)  \nonumber \\
& = &  o\left(1\right) + 
       \frac{1}{n} \sum_{s \in S} \sum_{j=1}^{N_{n-1}\left(s^{-1}\right)}   
                   \left(\ln \omega_{k_j\left(s^{-1}\right)}\left(s^{-1}\right) - \ln \omega_{k_j\left(s^{-1}\right)-1}\left(s\right)\right),
\end{eqnarray}
where for each $s \in S$, $k_1\left(s^{-1}\right), k_2\left(s^{-1}\right) \cdots, k_{N_{n-1}\left(s^{-1}\right)}\left(s^{-1}\right)$ are the 
time points $k$ when $Y_k = s^{-1}$ and
\begin{equation}
N_{n}\left(s\right) = \sum_{k=1}^{n} \bone\left(Y_{k} = s\right).
\end{equation}
Now consider the product space 
$\left(\Omega \times S^{\W}, \Bi_{\Omega} \otimes \Bi_{{\mathbb N}_0}, \P \otimes \mu\right)$. By
Theorems 6.5.5 and 6.6.1 of \cite{Dur10} we have $\P \otimes \mu$-a.s. for all $s \in S$,
\begin{equation}
\lim_{n \rightarrow \infty} \frac{N_{n}\left(s\right)}{n} = \frac{1}{d}.
\label{Equ:EFP}
\end{equation}
Further under assumption (A2) and using the Strong Law of Large Numbers for i.i.d. random variables we have 
$\P$-a.s., for every fixed $s \in S$, 
\[
\lim_{n \rightarrow \infty}    \frac{1}{N_{n-1}\left(s^{-1}\right)} 
\sum_{j=1}^{N_{n-1}\left(s^{-1}\right)} \ln {\omega_{k_j\left(s^{-1}\right)}\left(s^{-1}\right)}
= E\left[ \ln {\omega_{1}\left(s^{-1}\right)}\right],
\]
and also 
\[
\lim_{n \rightarrow \infty}    \frac{1}{N_{n-1}\left(s^{-1}\right)} 
\sum_{j=1}^{N_{n-1}\left(s^{-1}\right)} \ln {\omega_{k_j\left(s^{-1}\right)-1}\left(s\right)}
= E\left[ \ln {\omega_{1}\left(s\right)}\right].
\]
As $S$ is a symmetric set of generators for $G$, therefore $\P \otimes \mu$-a.s.,
\begin{equation}
\lim_{n \rightarrow \infty}\frac{1}{n} \ln \Phi (\eta_{n}) = \frac{1}{d}\sum_{s \in S}  
E\left[\ln {\omega_{1}\left(s^{-1}\right)} - \ln{\omega_{1}\left(s\right)}\right] = 0.
\label{Equ:Fundamental-Limit}
\end{equation}
So by Fubini's theorem, it follows that the equation~\eqref{Equ:Fundamental-Limit} holds
$\mu$-a.s., for every $\omega \in \Omega$ a.s. with respect to $P$.
Fix such an $\omega \in \Omega$. 
As $d \geq 3$ so find $\frac{1}{d-1} < \Delta < 1$.
Since almost sure convergence implies convergence in probability, so
$\exists \, M_{\mu}^{\omega} \in \Nbold$ such that for all $n \geq M_{\mu}^{\omega}$, 
\begin{equation}
\mu\left( \Phi\left(\eta_n\right) < \left(\frac{1}{\sqrt{\Delta}}\right)^{n}\right) > \frac{1}{2}.
\end{equation}
But recall that under $\mu$, the distribution of $\eta_n$ is uniform on the vertices of 
$\T^{n}_{d}$, so
\begin{equation}
\frac{\# \left\{\sigma_n \in \T_d \,\Big\vert\, 
     \Phi\left(\sigma_n\right) <  \left(\frac{1}{\sqrt{\Delta}}\right)^{n}\right\}}
     {d \left(d-1\right)^{n-1}}
> \frac{1}{2}
\end{equation}
for all $n \geq M_{\mu}^{\omega}$. 
%
%
Let $\beta_{n} = \Delta^{\frac{n}{2}}$. Observe that
$\displaystyle{\sum_{n=1}^{\infty} \beta^n < \infty}$. Now for $n \geq M_{\mu}^{\omega}$, 
\begin{equation}
     \sum_{\sigma_{n} \in \T^{n}_{d}}  \beta_{n} \left(\Phi\left(\sigma_{n}\right)\right)^{-1} 
\geq \mathop{\sum_{\sigma_{n}\in \T^{n}_{d}}}\limits_{\Phi\left(\sigma_n\right) <  \left(\frac{1}{\sqrt{\Delta}}\right)^{n}}  
     \beta_{n} \left(\Phi\left(\sigma_{n}\right)\right)^{-1}  
\geq \frac{1}{2} d \left(d-1\right)^{n-1}  \Delta^n.
\end{equation}
By the choice of $\Delta$ it follows that $P$-a.s. equation~\eqref{Equ:Main-Criterion} holds,
which completes the proof. 

\qed 

\subsection{Proof of Theorem~\ref{Thm:Speed}}
\label{SubSec:Proof-Speed}
Let $D_n := \left\vert X_n \right\vert$, then
\begin{eqnarray}
D_n & = & \sum_{i=1}^n \left( D_i - D_{i-1} \right) \nonumber \\
    & = & \sum_{i=1}^n \left(D_i - D_{i-1} - \Ebb_\omega^e\left[D_i - D_{i-1} \,\Big\vert\, X_0, \ldots , X_{i-1}\right]\right)  \nonumber \\
    &   & + \sum_{i=1}^n \Ebb_\omega^e\left[D_i - D_{i-1} \,\Big\vert\, X_0, \ldots , X_{i-1} \right]. \label{Equ:MG-Decom}
\end{eqnarray}
But then $\displaystyle{M_n := \sum_{i=1}^n \left(D_i - D_{i-1} - \Ebb_\omega^e\left[D_i - D_{i-1} \,\Big\vert\, X_0, \ldots , X_{i-1}\right]\right)}$
is a martingale with zero mean and bounded increments, so by Theorem 3 of \cite{Az67} 
\begin{equation}
\frac{M_n}{n} \rightarrow 0 \,\,\, \Pbb^e\mbox{-a.s..}
\label{Equ:MG-Conv}
\end{equation}  

Further it is easy to see that 
\[
D_i - D_{i-1} = \begin{cases}
                +1 & \mbox{if\ } X_{i-1} = e \\
                +1 & \mbox{if\ } X_{i-1} \not\in \left\{e, \overleftarrow{X}_{i-1} \right\} \\
                -1 & \mbox{if\ } X_{i-1} = \overleftarrow{X}_{i-1}.
                \end{cases}
\]
Thus,
\[
\Ebb_\omega^e\left[D_i - D_{i-1} \,\Big\vert\, X_0, \ldots , X_{i-1} \right] = 
1 - 2 \times \bone\left(X_{i-1} \neq e\right) \, \omega\left(X_{i-1}, \overleftarrow{X}_{i-1}\right).
\]

Now under our assumption (A3) with $\eps > \frac{1}{2\left(d-1\right)}$ there exists $\delta > 0$ such that
$P$-a.s
\[ 
\omega\left(x,\overleftarrow{x}\right) < \frac{1}{2} - \delta \left(d-1\right) \,\, \forall \,\, x \in \T_d.
\]
This is because 
$\displaystyle{\omega\left(x,\overleftarrow{x}\right) = 
1 -\sum_{x \sim y, \,\, y \neq \overleftarrow{x}} \omega\left(x,y\right)}$.
Thus $\Pbb^e$-a.s.
\begin{equation}
\liminf_{n \rightarrow \infty} \frac{1}{n} \sum_{i=1}^n \left(1 - 2 \times \bone\left(X_{i-1} \neq e\right) \,
                                                        \omega\left(X_{i-1}, \overleftarrow{X}_{i-1}\right) \right)  > 2 \delta \left(d-1\right) > 0.
\label{Equ:Rn_Conv}
\end{equation}
Finally, by \eqref{Equ:MG-Decom} 
$D_n = M_n +  
\displaystyle{\sum_{i=1}^n \left(1 - 2 \times \bone\left(X_{i-1} \neq e\right) \omega(X_{i-1}, \overleftarrow{X}_{i-1})\right)}$, so 
using equations \eqref{Equ:MG-Conv} and \eqref{Equ:Rn_Conv}
we conclude that \eqref{Equ:Speed} holds. 
\qed

\section*{Acknowledgment}
The authors would like to thank the two anonymous referees for their valuable and very detailed reviews which 
have improved the exposition of the article significantly. Thanks are also due to one of the referees for
pointing out the earlier work of Rozikov \cite{Rozi01} which was unknown to us.

\end{document}